\documentclass[12pt]{amsart}
\usepackage{amssymb,amsmath}
\newif\ifppddff
\ppddfftrue
\ifppddff\usepackage[pdftex]{graphicx}\else
\usepackage[dvips]{graphicx}\fi
\newcommand{\Graphics}[2]{\ifppddff \includegraphics[#1]{#2.pdf} \else \includegraphics[#1]{#2.eps}\fi}

      \newenvironment{changemargin}[2]{\begin{list}{}{
         \setlength{\topsep}{0pt}\setlength{\leftmargin}{0pt}
         \setlength{\rightmargin}{0pt}
         \setlength{\listparindent}{\parindent}
         \setlength{\itemindent}{\parindent}
         \setlength{\parsep}{0pt plus 1pt}
         \addtolength{\leftmargin}{#1}\addtolength{\rightmargin}{#2}
         }\item }{\end{list}}

\newcommand{\dof}[1]{{#1^\delta}}
\newcommand{\itdof}[2]{{#2^{\delta^{#1}}}}

\newcommand{\sfA}{\mathsf{A}}
\newcommand{\COR}{\mathsf{COR}}
\renewcommand{\ng}{\textsc{ng}}
\newcommand{\ns}{\textsc{ns}}
\newcommand{\wl}{\textsc{wl}}
\renewcommand{\sl}{\textsc{sl}}
\newcommand{\R}{\mathbb{R}}

\newcommand{\myfigure}[1]{\begin{figure}[!h]#1\end{figure}}

\newcommand{\lG}{\ell_\mathrm{G}}
\newcommand{\lR}{\ell_\mathrm{R}}

\newcommand{\lmin}{\ell}

\newcommand{\inv}{^{-1}}

\newcommand{\be}{\begin{enumerate}}
\newcommand{\ee}{\end{enumerate}}
\newcommand{\bi}{\begin{itemize}}
\newcommand{\ei}{\end{itemize}}
\newcommand{\itm}{\item}

\newtheorem{thm}{Theorem}
\newtheorem{prop}[thm]{Proposition}
\newtheorem{as}[thm]{Assumption}
\newtheorem{cor}[thm]{Corollary}
\newtheorem{lem}[thm]{Lemma}
\newtheorem{conj}[thm]{Conjecture}

\newtheorem{notat}[thm]{Notation}
\newtheorem{alg}[thm]{Algorithm}
\newtheorem{prob}[thm]{Problem}
\theoremstyle{definition}
\newtheorem{defn}[thm]{Definition}
\theoremstyle{remark}
\newtheorem{rem}[thm]{Remark}
\newtheorem{exa}[thm]{Example}

\newcommand{\bpf}{\begin{proof}}
\newcommand{\epf}{\end{proof}}

\newcommand{\blem}{\begin{lem}}
\newcommand{\elem}{\end{lem}}
\newcommand{\bthm}{\begin{thm}}
\newcommand{\ethm}{\end{thm}}
\newcommand{\bprp}{\begin{prop}}
\newcommand{\eprp}{\end{prop}}
\newcommand{\bcor}{\begin{cor}}
\newcommand{\ecor}{\end{cor}}
\newcommand{\bcon}{\begin{conj}}
\newcommand{\econ}{\end{conj}}
\newcommand{\bdfn}{\begin{defn}}
\newcommand{\edfn}{\end{defn}}
\newcommand{\bexm}{\begin{exa}}
\newcommand{\eexm}{\end{exa}}
\newcommand{\brem}{\begin{rem}}
\newcommand{\erem}{\end{rem}}
\newcommand{\bnot}{\begin{notat}}
\newcommand{\enot}{\end{notat}}
\newcommand{\balg}{\begin{alg}}
\newcommand{\ealg}{\end{alg}}

\renewcommand{\>}{\right > }

\long\def\forget#1\forgotten{} %

\begin{document}
\title[Random equations in Garside groups]{Solving random equations in Garside groups using length functions}

\author{Martin Hock}
\address{Department of Computer Science, University of Wisconsin, Madison, WI 53706, USA}
\email{mdhock@gmail.com}

\author{Boaz Tsaban}
\thanks{The second author was partially supported by the Koshland Center for Basic Research.}
\address{Department of Mathematics, Bar-Ilan University, Ramat-Gan 52900, Israel;
and
Department of Mathematics, Weizmann Institute of Science, Rehovot 76100, Israel}
\email{tsaban@math.biu.ac.il}
\urladdr{http://www.cs.biu.ac.il/\~{}tsaban}

\begin{abstract}
We give a systematic exposition of memory-length algorithms for
solving equations in noncommutative groups. This exposition
clarifies some points untouched in earlier expositions.
We then focus on the main ingredient in these attacks:
Length functions.

After a self-contained introduction to Garside groups,
we describe length functions induced by the
greedy normal form and by the rational normal form in these
groups, and compare their worst-case performances.

Our main concern is Artin's braid groups, with their
two known Garside presentations, due to
Artin and due to Birman-Ko-Lee (BKL).
We show that in $B_3$ equipped with the BKL presentation,
the (efficiently computable) rational normal form of each element
is a geodesic, i.e., is a representative of minimal length for that element.
(For Artin's presentation of $B_3$, Berger supplied in 1994 a method to obtain geodesic
representatives in $B_3$.)

For arbitrary $B_N$, finding the geodesic length of an element is NP-hard, by a 1991 result of by Paterson and Razborov.
We show that a good estimation of the geodesic length of an element of $B_N$
in Artin's presentation is measuring the length of its rational form in the \emph{BKL} presentation.
This is proved theoretically for the worst case, and experimental evidence
is provided for the generic case.
\end{abstract}

\maketitle

\section{Solving random equations}

All groups considered in this paper are multiplicative noncommutative groups, with an efficiently solvable word problem,
that is, there is an efficient algorithm for deciding whether two given (finite products of) elements in the group are equal
as elements of the group. Throughout this paper, $G$ denotes such a group.

Problems involving solutions of equations in groups have a long history,
and are nowadays also explored towards applications
in public-key cryptography \cite{GBC}. We mention some of the more elegant problems of this type.

\begin{prob}[Conjugacy Search]\label{CS}
Given conjugate $a,b\in G$, find $x\in G$ such that $b=xax\inv$.
\end{prob}

\begin{prob}[Root Search]\label{RS}
Given $a\in G$, find $x\in G$ such that $a=x^2$, provided
that such $x$ exists.
\end{prob}

\begin{prob}[Decomposition Search]\label{DS}
Let $H$ be a proper subgroup of $G$.
Given $a,b\in G$, find $x,y\in H$ such that $b=xay$, provided
that there exist such $x,y$.
\end{prob}

We will discuss the meaning of the terms ``given'' and ``find'', appearing in Problems \ref{CS}--\ref{DS}, later.

Problems \ref{CS}--\ref{DS}, as well as many additional ones, can be stated generally as follows.
By \emph{free-group word} $w(t_1,\dots,t_k)$ we mean a product of \emph{variables}
$t_{i_1}^{\epsilon_1}\cdot t_{i_2}^{\epsilon_2}\cdot\ldots\cdot t_{i_n}^{\epsilon_n}$
for any choice of a positive integer $n$ and elements $i_1,\dots,i_n\in\{1,\dots,k\}$ and $\epsilon_1,\dots,\epsilon_n\in\{1,-1\}$,
such that no cancellation is possible, that is, for each $j=1,\dots,n$, if $i_j=i_{j+1}$, then $\epsilon_j\neq -\epsilon_{j+1}$.

\begin{prob}[Solution Search]\label{SS}
Fix $H_1,\dots,H_k\le G$ and a free-group word $w(t_1,\dots,t_{k+n})$.
Given parameters $p_{1},\dots,p_n\in G$ and an element $c\in G$, find $x_1\in H_1,\dots,x_k\in H_k$ such that
$c=w(x_1,\dots,x_k,\allowbreak p_{1},\dots,p_n)$, provided that there exist such $x_1,\dots,\allowbreak x_k$.
\end{prob}

Problem \ref{SS} deals with the solution of a single solvable equation (with parameters).
It can also be stated for systems of \emph{several equations}.
The algorithms proposed here easily generalize to cover this case, cf.\ \cite{BraidEqns}.

\subsection{Making the problems meaningful}
It suffices to discuss Problem \ref{SS}.

First, all given information must be coded in some compact form.
For example, the subgroups $H_1,\dots,H_k$ of $G$ may be described by lists of generators and
relations, all (the list, the generators, and the relations) of manageable length.

Second, the problem may require that it be possible to find a solution for each possible instance of
the problem, or for a certain portion of the instances.
Already in the case of free groups, the problem of solving equations in this sense is extremely difficult.
For example, the problem of solving quadratic equations over free groups is known to be NP-hard.

Alternatively, the instances of the problem may be chosen according to a certain distribution $D$, and
we may require that a solution can be found with a high-enough probability (a \emph{probabilistic model}).

Finally, by ``find'' we mean ``find efficiently'', i.e., use an algorithm with a feasible running time.
Otherwise, in most cases of interest the problems are solvable. E.g., if $G$ is a finitely generated group
with solvable word problem, then we can solve Problem \ref{SS} by enumerating $G^k$ recursively, and
trying all possible solutions until one is found. This algorithm always succeeds in a finite running time,
but usually this running time is infeasible.

In this discussion, all quantitative terms (compact, efficient, significant, etc.) have two natural interpretations:
Concrete (e.g., of size less than 1GB) or asymptotic (e.g., polynomial in the size of the input).

\subsection{The probabilistic model}
With an eye towards applications, we will always use the probabilistic version of the problems,
where we wish to find (efficiently) a solution with a significant probability,
provided that the instances of the problem are chosen according to a certain known distribution $D$.

More precisely, in Problem \ref{SS} we fix a distribution $D$ on $G^{k+n}$ such that for each $(x_1,\dots,x_k,p_{1},\dots,p_n)$ in
the support of $D$, we have that $x_1\in H_1,\dots,x_k\in H_k$. An instance of the problem is generated as follows:
A secret tuple $(x_1,\dots,x_k,p_{1},\dots,p_n)\in G^{k+n}$ is chosen according to the distribution $D$, and
we are given $p_{1},\dots,p_n$ and an element $c\in G$ equal to $w(x_1,\dots,x_k,p_{1},\dots,p_n)$ in $G$.
We must then search for elements $\tilde x_1\in H_1,\dots,\tilde x_k\in H_k$ such that with a significant probability,
$c=w(\tilde x_1,\dots,\tilde x_k,p_{1},\dots,p_n)$ in $G$.

By peeling off known parameters on the left of the given word $w(x_1,\allowbreak\dots,x_k, p_{1},\dots,p_n)$,
we may assume that it begins with a variable $x_i$ (possibly inverted).
If we are able to find $x_i$ (with a significant probability), we can treat it as a parameter henceforth,
and proceed to the next leading variable after peeling off all parameters on the left.
Continuing in this manner, we find suggestions for all variables, and can check whether we
obtained a solution.

Thus, it is natural to consider the following problem.

\begin{prob}[Leading-Variable Search]\label{LVS}
Fix $H_1,\dots,H_k\le G$ and a free-group word $t_1\cdot w(t_1,\dots,t_{k+n})$.
Given parameters $p_{1},\dots,p_n\in G$ and an element $c=x_1\cdot w(x_1,\dots,x_k,\allowbreak p_{1},\dots,p_n)\in G$
such that $x_1\in H_1,\dots,x_k\in H_k$, find $\tilde x_1\in H_1$, such that there are
$\tilde x_2\in H_2,\dots,\tilde x_k\in H_k$ with
$c=\tilde x_1\cdot w(\tilde x_1,\dots,\tilde x_k,\allowbreak p_{1},\dots,p_n)$.
\end{prob}

Clearly, any algorithm solving Problem \ref{SS} also solves Problem \ref{LVS}, with at least the same
probability of success. On the other hand, an algorithm for Problem \ref{LVS} can be iterated, as explained
above, to obtain a solution for Problem \ref{SS} (with a smaller probability of success, which also depends on its
performance on the induced distributions along the iteration).

\subsection{Decision problems}
All mentioned problems also have a \emph{decision} version.
For example, the \emph{Congugacy Problem} is: Given $a,b\in G$, are they conjugate?
If we only consider algorithms with \emph{bounded running time},
then a solution to the search version also implies a solution to the decision version, in the following sense.

Assume that $\sfA$ is an algorithm searching for solutions of equations of a certain type (e.g., $b=xax\inv$),
and that its running time is bounded, say by a certain function of the length of its input.
We define a \emph{decision} algorithm $\sfA'$ with running time bounded by the same function:
Given an instance of the equation to be checked, run $\sfA$ on this instance until the running time reaches its bound,
and then terminate it if it did not terminate already.
If a solution was found, the decision of $\sfA'$ is \emph{Yes}. Otherwise, it is \emph{No}.

Assume that the instances of the equation are distributed according to some
distribution $E$. This induces a distribution $D$ on the \emph{solvable} equations,
by conditioning that the chosen equation be solvable.
Let $p$ be the probability that $\sfA$ finds a solution to (necessarily, solvable) equations distributed according to $D$.

For each specific instance of the equation, $\sfA'$ is correct in probability at least $p$:
If this instance has a solution, it will be found by $\sfA$ in probability $p$, in which case $\sfA'$ decides ``Yes''.
And if this instance has no solution, then in probability $1$, $\sfA$ will not find a solution (because there is none),
and $\sfA'$ decides ``No''.

This can also be viewed as follows: Let $q=1-p$.
The probability that $\sfA'$ comes up with a wrong answer is:
\begin{eqnarray*}
\lefteqn{P(\mbox{Wrong decision})=}\\
& = & P(\mbox{Decision}=\mbox{Yes} \mid \nexists\mbox{Solution})\cdot P(\nexists\mbox{Solution})+\\
& + & P(\mbox{Decision}=\mbox{No} \mid \exists\mbox{Solution})\cdot P(\exists\mbox{Solution}) =\\
& = & 0\cdot P(\nexists\mbox{Solution}) + q \cdot P(\exists\mbox{Solution}) =\\
& = & q \cdot P(\exists\mbox{Solution}).
\end{eqnarray*}

In particular, this probability is at most $q$, and the worst case is when $P(\exists\mbox{Solution})$ is $1$, in which the
distribution may be assumed to be supported by solvable instances, and we are actually in the \emph{search} version of the problem.

This justifies, to some extent, restricting attention to search problems when working in the probabilistic model, with
algorithms of bounded running time.

\section{The memory-length approach}

The potential usefulness of length functions for solving the conjugacy search problem was identified in \cite{HT}.
In \cite{LBCS, BraidEqns}, it was pointed out that this approach can be used to solve arbitrary (systems of) equations.

Let $H\le G$ be generated by elements $a_1,\dots,a_m$ of $G$.
Assume that an instance $x\cdot w(x,x_2,\dots,x_k,\allowbreak p_{1},\dots,p_n)$ of Problem \ref{LVS} is
chosen according to a certain distribution $D$, with $H_1=H$,
and we are given $c$ which is equal to it in $G$.
Let $w=w(x,x_2,\dots,x_k,\allowbreak p_{1},\dots,p_n)$.

Let $A=\{a_1,\dots,a_m\}^{\pm 1}$.
Assume that the shortest expression of $x$ as a product of elements of $A$ has length $n$.
Let $\COR(x)$ be the set of all $a\in A$ which appear first in an expression of $x$ as a product of $n$ generators,
i.e., $\{a\in A : x\in_G a\cdot A^{n-1}\}$.
For each $a\in\COR(x)$, $a\inv x$ has an expression of length $n-1$, whereas for $a\notin\COR(x)$,
$a\inv x$ may in general not have an expression shorter than $n+1$. In particular, we expect
$a\inv x$ to be ``shorter'' when $a\in\COR(x)$ than when $a\notin\COR(x)$.
Heuristically, this expectation is extended to $xw$.

Often, we cannot compute the length of a shortest expression of a group element, and we only
assume that we have an efficiently computable function $\ell:G\to\R_{\ge 0}$, which approximates
the above situation, i.e., such that
$\ell(abw)$ tends to be greater than $\ell(w)$ for $w\in G,a,b\in\{a_1,\dots,a_m\}^{\pm 1}$.

By standard arguments, we may for convenience assume that $n$ is known \cite{BraidEqns, Thomp1}.\footnote{This
has a computational cost, so we cannot assume that we know the lengths of shortest expressions of many elements.}
One may then try all $a\in A$, and pick one with $\ell(a\inv xw)$ minimal. Hopefully,
$a\in\COR(x)$, and we can continue with the peeled-off element $a\inv xw$. After $n$ steps, we hopefully have
(a shortest expression for) $x$.

In cases of interest this approach does not work as stated \cite{LBCS}, and the following
improvement was proposed in \cite{BraidEqns}.

\subsection{The memory-length algorithm}
Using the above-mentioned notation, the algorithm generates an ordered list of $M$ sequences
of length $n$, with the aim that with a significant probability, a sequence
$$((j_1,\epsilon_1),(j_2,\epsilon_2), \dots, (j_n,\epsilon_n)),$$
such that $x=a_{j_1}^{\epsilon_1}a_{j_2}^{\epsilon_2}\dots a_{j_n}^{\epsilon_n}$ in $G$,
appears in the list, and tends to be among its \emph{first} few members. It consists of the following steps:

\subsubsection*{Step 1}
For each $j=1,\dots,m$ and each $\epsilon\in\{1,-1\}$, compute
$a_j^{-\epsilon}c = a_j^{-\epsilon}xy$, and give $(j,\epsilon)$
the score $\ell(a_j^{-\epsilon}c)$. Keep in
memory the $M$ elements $(j,\epsilon)$ with the best (=lowest) scores.

\subsubsection*{Steps $s=2,3,\dots,n$}
For each sequence $((j_1,\epsilon_1), \dots, (j_{s-1},\epsilon_{s-1}))$
out of the $M$ sequences stored in the
memory, each $j_s=1,\dots,m$, and each $\epsilon_s\in\{1,-1\}$,
compute
$$\ell(a_{j_{s}}^{-\epsilon_{s}}(a_{j_{s-1}}^{-\epsilon_{s-1}}\cdots a_{j_1}^{-\epsilon_1}c)) =
\ell(a_{j_{s}}^{-\epsilon_{s}}a_{j_{s-1}}^{-\epsilon_{s-1}}\cdots a_{j_1}^{-\epsilon_1}xy),$$
and assign this score to the sequence
$((j_1,\epsilon_1), \dots,(j_{s},\epsilon_{s}))$.
Keep in memory only the $M$ sequences with the best scores.

The algorithm terminates after $n$ steps, with $M$ proposals for
$((j_1,\epsilon_1),\allowbreak(j_2,\epsilon_2), \dots, (j_n,\epsilon_n))$.

\medskip

It is not difficult to see that the complexity of this algorithm is $n(n+4m+1)M/2$ group operations and evaluations of $\ell$.

It is interesting to note that this algorithm may also be useful for solving the following.

\begin{prob}[(Shortest) Subgroup Membership Search]
Given $a_1,\dots,\allowbreak a_m\in G$ and $x\in \<a_1,\dots,a_m\>$, find a (shortest possible)
expression of $x$ as a product of elements from the set $\{a_1,\dots,a_m\}^{\pm 1}$.
\end{prob}

\subsection{Sufficiency for the general problem}
Assume that the algorithm succeeds, with a significant probability,
to have the leading element $x$ in the final list.
Then we have the following.

If there is only \emph{one unknown variable} in the equation (e.g., Problems \ref{CS}--\ref{DS}), then
we can check (in running time $M$) all elements in the list and find one which is a solution to the problem.

In the general case (Problem \ref{SS}) there are several unknown variables,
and we can iterate the algorithm by checking each suggestion in the list.
The overall complexity is in principle $M^k$.
However, the suggestions for each variable are ordered more or less according to their likelihood,
and it suffices to check, for some $N\ll M$, the $N$ most likely solutions.
This reduces the complexity to $N^k$, or more precisely to $N_1\cdot N_2\cdots N_k$, where $N_k$
is the number of elements required at the $k$th step, and it is likely that $N_{i+1}\ll N_i$ for each $i$.

\subsection{Improvements}
Certain simple modifications in the memory-length algorithm increase its success rates.
We refer the reader to \cite{Thomp1} for details.

\subsection{The length function}
For this algorithm to be meaningful and useful, one must have a
good and efficiently computable length function on the group $G$.
Our introduction of the memory-length algorithm suggests a natural
model for comparing length functions for appropriateness to this method.
We explore this below, after introducing a
new proposal for a length function on the braid group.
The braid group is, thus far, the most popular in applications related to cryptography \cite{GBC}.
Most of these cryptographic applications give rise to an equation, whose solution
would imply the insecurity of the application. Thus, it is natural to look for good length
functions on this group. See \cite{GBC} for more details.

\section{Excursion: Garside groups}

We are going to consider two Garside structures on the braid group
(to be defined). This section is an essentially self-contained introduction
to Garside groups, and may be skipped by readers who are familiar with
this concept, and by readers who do not insist on understanding
all details of this paper.

Garside groups were introduced by Dehornoy and
Paris \cite{DP}, and later in a more general form by Dehornoy
\cite{Deh02}. We treat the latter, more general case. All unproved
assertions, as well as most of the proved ones, are from \cite{DP}.

\subsection{Garside Monoids and Groups}
Let $M$ be a monoid with cancellation.
$x\in M$ is an \emph{atom} if $x\ne 1$, and $x=ab$ for $a,b\in M$ implies $a=1$ or $b=1$.
$M$ is \emph{atomic} if $M$ is generated by its atoms, and for each $a\in M$,
the maximum number of atoms in an expression of $a$ as a product
of atoms, denoted $\|a\|$, exists.
It follows that $\|ab\|\ge\|a\|+\|b\|$ for all $a,b\in M$.
In particular, as $1=1\cdot 1$, we have that $\|1\|\ge \|1\|+\|1\|$, and thus $\|1\|=0$.
For $a\neq 1$, $\|a\|>0$.

Let $M$ be an atomic monoid. For $a,b\in M$, $a$ is a \emph{left divisor} of $b$
if there is $c\in M$ such that $ac=b$.
Similarly, $a$ is a \emph{right divisor} of $b$
if there is $c\in M$ such that $ca=b$.
$a\in M$ is a \emph{Garside element} of $M$ if
its left divisors and right divisors coincide,
and include all atoms of $M$.

$M$ is a \emph{Garside monoid} if it is atomic, has a Garside
element, and for all $a,b\in M$, a greatest common divisor $a\wedge b$
and a least common multiple $a\vee b$ of $a$ and $b$ exist in $M$, both with
respect to left divisibility.

For $a,b\in M$, the \emph{complement} $a\setminus b$ is the
unique $c\in M$ such that $ac=a\vee b$.
The closure of the set of atoms under the operations of
complement and least common multiple is the set $S$ of \emph{simple elements} of $M$.
The least common multiple of all elements of $S$, if it exists (e.g., if $M$ is finitely generated),
is called the \emph{fundamental element} of $M$ and denoted $\delta$.
$\delta$, if it exists, is the least Garside element of $M$.

$G$ is a \emph{Garside group} if it is the group of fractions of a Garside monoid $M$.
In this case, the elements of $M$ are called the \emph{positive elements} of $G$.
In the remainder of this section, $M$ is a Garside group with a fundamental element $\delta$, and $G$ is the Garside group of fractions of $M$.

\subsection{Greedy Normal Form}
For $x\in M$ with $x\neq 1$, the simple element $\delta\wedge x\neq 1$.
Define $\partial(x)=(\delta\wedge x)\inv x$. Then $\partial(x)\in M$,
and as $x = (\delta\wedge x)\partial(x)$, $\|x\| \ge \|\delta\wedge x\|+\|\partial(x)\|>\|\partial(x)\|$.
Define simple elements $s_1,s_2,\dots$, as follows. Set $x_1=x$, and
for each $i=1,\dots,r$, let $s_i = \delta \wedge x_{i}$, and $x_{i+1}=\partial(x_{i})$.
$\|x\|=\|x_1\|>\|x_2\|>\dots\ge 0$, and thus there is a minimal $n$ such that $x_{n+1}=1$.
$x = s_1\cdots s_n$.
Let $k\ge 0$ be maximal with $s_i=\delta$, and define $p_i=s_{k+i}$, $i=1,..,r$, $r=n-k$.
The expression
$$x = \delta^k p_1\cdots p_r$$
is called the \emph{greedy normal form} of $x$.

Consider now $x\in G\setminus M$. If $x=\delta^k s$ and $s\in M$, then $k<0$.
Take the maximal integer $k$ such that $x=\delta^k s$ for some $s\in M$.
Fix such $s$, and let $\delta^0p_1\cdots p_r=p_1\cdots p_r$ be the greedy
normal form of $s$.
The greedy normal form of $x$ is then again defined to be $\delta^k p_1\cdots p_r$.

By the construction, we have that
$p_{i+1} \wedge p_i^{-1}\delta = (p_{i+1}\cdots p_r \wedge \delta) \wedge p_i^{-1}\delta
= p_{i+1}\cdots p_r \wedge (\delta \wedge p_i^{-1}\delta) =
x_{i+1} \wedge p_i^{-1}\delta = 1$ for all $i=1,\dots,r-1$, and that
$p_r\neq 1$. We say in such cases that the sequence $p_1,\dots,p_r$ is \emph{left-weighted}.

\subsection{Rational Normal Form}
Following Thurston \cite[Chapter 9]{Eps}, Dehornoy and Paris define the
\emph{rational normal form}\footnote{Also called \emph{mixed} or \emph{symmetric} normal form.}
of an element $x\in G$. To this end, we need the following.

\bthm[Dehornoy-Paris \cite{DP}]\label{quotient}
For each $x\in G$, there is a unique pair $(u,v)$ in $M\times M$
such that $x=u\inv v$ and $u\wedge v=1$.
\ethm


Let $x\in G$, and let $u,v\in M$ be as in Theorem \ref{quotient}.
Let $s_1\cdots s_k$, $p_1\cdots p_l$ be the greedy normal form
of $u,v$, respectively. The rational normal form of $x$
is the expression
$$x=(s_1\cdots s_k)\inv (p_1\cdots p_l).$$
All $s_i,p_j$ are simple, $s_1\wedge p_1=1$, and
the sequences $s_1,\dots,s_k$ and $p_1,\dots,p_l$ are both
left-weighted. (The special cases where $k=0$ or $l=0$ are also allowed.)

For each $a\in G$, define $\tau(a)=\dof{a}=\delta\inv a\delta$.
$\tau$ is an inner automorphism of $G$, and its
$n$th iterate at $a$ is $\tau^n(a)=\itdof{n}{a}$.
$\tau$ maps simple elements to simple elements:
For each simple $s$, let $p$ be such that $sp=\delta$. Then $p$ is simple,
and thus there is a simple $q$ with $pq=\delta$. Then
$$s\delta=spq=\delta q,$$
and thus $\dof{s}=q$ is simple.
In particular, $M$ is invariant under $\tau$.
Any automorphism of $G$ mapping positive elements to positive elements,
maps atoms to atoms.
It follows that $\tau$ is a permutation of the atoms of $M$.

One can obtain the rational normal form from the greedy normal form.
To see this, we use the following.

\blem\label{lwlem}
If $s,p$ are simple and $sp$ is left-weighted, then so are $\dof{s}\dof{p}$ and $\itdof{-1}{s}\itdof{-1}{p}$.
\elem
\bpf
If $ac=b$ are all positive, then
$\itdof{\pm 1}{a}\itdof{\pm 1}{c}=\itdof{\pm 1}{(ac)}=\itdof{\pm 1}{b}$, and $\itdof{\pm 1}{c}\in M$.
Thus, $\tau^{\pm 1}$ both map left divisors to left divisors, and therefore
$$\itdof{\pm 1}{(a\wedge b)}=\itdof{\pm 1}{a}\wedge \itdof{\pm 1}{b}$$
for all $a,b\in M$.
Now, assume that $sp$ is left-weighted. Then
$$(\itdof{\pm 1}{s})\inv\delta\wedge\itdof{\pm 1}{p} =
\itdof{\pm 1}{(s\inv\delta)}\wedge\itdof{\pm 1}{p} =
\itdof{\pm 1}{(s\inv\delta\wedge p)} =
\itdof{\pm 1}1=1,$$
showing that $\itdof{\pm 1}{s}\itdof{\pm 1}{p}$ is left-weighted.
\epf

\bprp\label{ver}
If $s,p$ are simple and $sp$ is left-weighted, then so are $((\itdof{k}{p})\inv\delta)((\itdof{k+1}{s})\inv\delta)$,
for all integer $k$.
\eprp
\bpf
Assume that $sp$ is left-weighted.
Then so is $(p\inv\delta)((\dof{s})\inv\delta)$:
$$(p\inv\delta)\inv\delta\wedge((\dof{s})\inv\delta)=\dof{p}\wedge\dof{(s\inv\delta)} =
\dof{(p\wedge(s\inv\delta))}=\dof1=1.$$
By Lemma \ref{lwlem},
$((\itdof{k}{p})\inv\delta)((\itdof{k+1}{s})\inv\delta)
= \itdof{k}{((p\inv\delta)((\dof{s})\inv\delta))}$
is also left-weighted.
\epf

Let $\delta^k p_1\cdots p_r$ be the greedy normal form of $x$.
Consider three possible cases.

\subsubsection*{Case 1: $k\ge 0$}
Then $\delta^k p_1\cdots p_r$ is already a rational normal form (with a trivial negative part).

\subsubsection*{Case 2: $k=-m<0$ and $m\ge r$}
By definition, $\delta^{-n}a = \itdof{n}a\delta^{-n}$ for all $a$ and all $n$.
Using this, we have that
\begin{eqnarray*}
\lefteqn{\delta^{-m} p_1\cdots p_r = \delta\inv\itdof{m-1}{p_1}\cdot\delta\inv\itdof{m-2}{p_2}\cdot\ldots\cdot \delta\inv\itdof{m-r}{p_{r}}\cdot\delta^{-(m-r)}=}\\
& = & \left(\delta^{m-r}\cdot(\itdof{m-r}{p_{r}})\inv\delta \cdot\ldots\cdot(\itdof{m-2}{p_2})\inv\delta \cdot (\itdof{m-1}{p_1})\inv\delta\right)\inv.
\end{eqnarray*}
By Proposition \ref{ver}, the last inverted expression is left-weighted, and thus
we have a rational form, with a trivial positive part.

\subsubsection*{Case 3: $k=-m<0$ and $m<r$}
In the same manner, we have that
\begin{eqnarray*}
\lefteqn{\delta^{-m} p_1\cdots p_r = \delta\inv\itdof{m-1}{p_1}\cdot\delta\inv\itdof{m-2}{p_2}\cdot\ldots\cdot \delta\inv p_{m} \cdot p_{m+1}\cdot\ldots\cdot p_r=}\\
& = & \left(p_{m}\inv\delta \cdot\ldots\cdot(\itdof{m-2}{p_2})\inv\delta \cdot (\itdof{m-1}{p_1})\inv\delta\right)\inv (p_{m+1}\cdot\ldots\cdot p_r),
\end{eqnarray*}
By Proposition \ref{ver}, each of the bracketed expressions is left-weighted.
Thus, this expression is in rational normal form.

\section{Several length functions on Garside groups}

Let $M$ be a Garside monoid with fundamental element $\delta$, and $G$ be its
group of quotients.

\begin{as}\label{asm}
We assume that for each simple $s\in M$, the minimal length $\lmin(s)$ of an expression
of $s$ as a product of atoms can be efficiently computed.
\end{as}

There is always an algorithm for computing $\lmin(s)$: Enumerate all words of
length $1,2,3,\dots$, until one equal to $s$ is found. The running time is
bounded by $k^{\lmin(a)}\le k^{\|a\|}$, where $k$ is the number of atoms. But this is in general
infeasible. When Assumption \ref{asm} fails, one may use in applications an estimation of $\lmin$ instead
of the true function.

Fortunately, in the specific monoids in which we are interested,
all relations are length-preserving, and thus $\lmin(s)$ is just
the length of any expression of $s$ as a product of atoms.
Thus, Assumption \ref{asm} is true in our applications.

\bexm[Artin's presentation of $B_N$]
\label{ArtinBN}
Consider the monoid $B_N^+$ generated by $\sigma_1,\dots,\sigma_{N-1}$,
subject to the relations
\begin{eqnarray*}
\sigma_i \sigma_{i+1}\sigma_i & = & \sigma_{i+1} \sigma_i \sigma_{i+1};\\
\sigma_i \sigma_j          & = & \sigma_j \sigma_i \mbox{ when }|i-j|>1.
\end{eqnarray*}
The quotient group of this monoid is the braid group $B_N$ on $N$ strings.
$B_N^+$ is a Garside monoid with atoms $\sigma_1,\dots,\sigma_{N-1}$,
and fundamental element
$$\delta=(\sigma_1\cdots\sigma_{N-1})(\sigma_1\cdots\sigma_{N-2})\cdots(\sigma_1\sigma_2)\sigma_1.$$
The positive elements of $B_N$ are the words in
$\sigma_1,\dots,\sigma_{N-1}$ not involving inverses of generators.
As the relations are length preserving, all expressions of a positive element as a product of atoms
have the same length. Thus, for $a\in M$, $\|a\|$ is the length of a (any) presentation of $a$.

Elements of $B_N$ can be identified with braids having $N$ strings,
where each generator $\sigma_i$ performs a half-twist on the $i$th and $i+1$st strings.
This way, $\delta$ is a half-twist of the full set of strings.
The simple elements correspond to positive braids in which
any two strings cross at most once.  A simple element is described
uniquely by the permutation it induces on the strings, and every permutation
of the $N$ strings corresponds to a simple element.
\eexm

\begin{exa}[BKL presentation of $B_N$]
\label{BKLBN}
Generalizing the geometric interpretation in Example \ref{ArtinBN}
to allow half-twists of the $i$th and the $j$th string for arbitrary
$i,j$, Birman, Ko, and Lee \cite{BKL98} introduced the following presentation
of the braid group $B_N$.
The monoid $BKL_N^+$ is generated by $a_{t,s}$, $1\le s<t\le N$,
subject to the relations
$$     \begin{array}{l}
       a_{t,s} a_{r,q} = a_{r,q} a_{t,s}
                     \quad\mbox{if}\quad (t-r)(t-q)(s-r)(s-q) > 0; \\
       a_{t,s} a_{s,r} = a_{t,r} a_{t,s} = a_{s,r} a_{t,r}
                     \quad\mbox{if}\quad t > s > r.
     \end{array}
$$
Also here, the relations are length preserving, and thus the norm is equal to
the number of atoms in any expression of the element.

This monoid also has the braid group $B_N$ as its quotient group.
In terms of Artin's presentation (Example \ref{ArtinBN}),
the Birman-Ko-Lee (BKL) generators can be expressed by
$$a_{t,s} = (\sigma_{t-1}\cdots\sigma_{s+1}) \sigma_s
 (\sigma_{s+1}^{-1}\cdots\sigma_{t-1}^{-1}).$$
$BKL_n^+$ is a
Garside monoid with fundamental element
$$\delta=a_{n,n-1}a_{n-1,n-2}\cdots a_{2,1}.$$
Here too, a simple element is described uniquely by the permutation it
induces on the strings. However, not every permutation of the $n$ strings
corresponds to a simple element.
\end{exa}

\bdfn
Let $M$ be a Garside monoid with Garside group $G$, and let $x\in G$.
\be
\itm $\lmin(x)$, the \emph{minimal length} of $x$, is the minimal length of an expression
of $x$ as a product of elements of $A^{\pm1}$, where $A$ is the set of atoms of $M$.
\itm $\lG(x)$, the \emph{greedy length} of an $x$,
is the sum of the minimal lengths of all simple elements (including the inverted ones)
in the greedy normal form of $x$. Similarly:
\itm $\lR(x)$, the \emph{rational length} of $x$,
is the sum of the minimal lengths of all simple elements (including the inverted ones)
in the rational normal form of $x$.
\ee
\edfn

Specifically,
if the greedy normal form of $x$ is $\delta^k s_1\cdots s_r$,
then $\lG(x)=k\cdot\lmin(\delta)+\lmin(s_1)+\dots+\lmin(s_r)$,
and if the rational normal form of length of $x$ is
$(s_1\dots s_k)\inv p_1\dots p_l$, then $\lR(x)=\lmin(s_1)+\dots+\lmin(s_k)+\lmin(p_1)+\dots+\lmin(p_l)$.

\bprp\label{samelen}
For each $a\in M$, $\lmin(\dof{a})=\lmin(a)$.
\eprp
\bpf
Let $n=\lmin(a)$, and $a=a_1\cdots a_n$ with $a_1,\dots,a_n$ atoms.
Then $\dof{a}=\dof{a_1}\cdots \dof{a_n}$. As conjugation by $\delta$ moves
atoms to atoms, $\lmin(\dof{a})\le n = \lmin(a)$.
Similarly, if $m=\lmin(\dof{a})$ and $\dof{a}=b_1\cdots b_m$ with $b_1,\dots,b_m$ atoms,
then $a = \itdof{-1}{\dof{a}}=\itdof{-1}{b_1}\cdots \itdof{-1}{b_m}$, and
as conjugation by $\delta$ moves atoms to atoms, $\lmin(a)\le m = \lmin(\dof{a})$.
\epf

The presentation in the previous section of the rational normal form in terms of
the greedy normal form gives the following.
\bcor
The rational length of an element with greedy normal form $\delta^{-m} s_1\cdots s_r$, where $0<m\le r$,
is
$$\lmin(s_1\inv\delta)+\dots+\lmin(s_m\inv\delta)+\lmin(s_{m+1})+\dots+\lmin(s_r),$$
and similarly for the cases where $m\le 0$ or $0<r<m$.
\ecor

\bcor\label{removing}
If the relations of $M$ are length-preserving, then
the rational length of an element with greedy normal form $\delta^k s_1\cdots s_r$ can be obtained
by removing $2 \sum_{i=1}^{\min(r,k)}{\lmin(s_i)}$ from its greedy normal length.
\ecor
\bpf
If the relations of $M$ are length-preserving, we have that $\lmin(ab)=\lmin(a)+\lmin(b)$ for all
$a,b\in M$, and thus for simple $s$, $\lmin(\delta)=\lmin(s)+\lmin(s\inv\delta)$, that is,
$\lmin(s\inv\delta)=\lmin(\delta)-\lmin(s)$.
\epf

This shows, in particular, that the length function considered in \cite{LBCS, BraidEqns}
in the case of the Artin presentation of $B_N$ is in fact the rational length for the Artin presentation
of $B_N$. This was first pointed out to us by Dehornoy.

\subsection{Quasi-geodesics in Garside groups}
Even when the relations are length-preserving,
it is generally not the case that an efficient algorithm for computing the minimal length $\lmin(x)$
is available. Even if the monoid relations are length-preserving, finding $\lmin(x)$ for
$x$ not in the monoid (nor in its inverse) may be a difficult task.
Indeed, assuming $P\neq NP$, there is no polynomial-time algorithm computing $\lmin(x)$ with respect to
the Artin presentation of $B_N$, for arbitrary $N$ and $x\in B_N$ \cite{PatRaz}.
Fortunately, in Garside groups $\lmin(x)$ can be \emph{approximated}. For simplicity, we treat
the case of length-preserving relations, so that $\lmin$ is easy to compute on positive elements.

\begin{thm}\label{basiclengths}
Let $M$ be a Garside monoid with length preserving relations and fundamental element $\delta$, and
let $G$ be its fractions group. For each $x\in G$:
\be
\itm If $x\in M$, then $\lG(x)=\lR(x)=\lmin(x)$.
\itm If $x\in M\inv$, then $\lR(x)=\lmin(x)$.
\itm $\lmin(x)\le \lR(x)\le \lG(x)\le (2\lmin(\delta)-1)\lmin(x)$.
\itm $\lR(x)\le (\lmin(\delta)-1)\lmin(x)$.
\ee
Moreover, these bounds in (3) cannot be improved.
\end{thm}
\begin{proof}
(1) For $x\in M$, each normal form gives some positive presentation of $x$, and thus the corresponding
length is the same as the minimal length.

(2) Fix $x\in M^{-1}$. Then $\lR(x)=\lR(x\inv)$, and by (1), $\lR(x\inv)=\lmin(x\inv)=\lmin(x)$.

(3) The first inequality is clear. The second follows from Corollary \ref{removing}.
We prove the third.
Let
\begin{equation}\label{minpres}
x=a_1^{\epsilon_1}\cdots a_m^{\epsilon_m}
\end{equation}
with $m=\lmin(x)$, $a_1,\dots,a_m$ atoms, and $\epsilon_1,\dots,\epsilon_m\in\{1,-1\}$.
For each atom $a$, let $\bar a$ be the simple element such that $\bar a a=\delta$.
Then $a\inv = \delta\inv\bar a$. Rewrite each negative atom in the equation \ref{minpres}
in this form, and move all occurrences of $\delta\inv$ to the left, using the relation
$a\delta\inv  = \delta\inv a^{\delta\inv}$.
Let $n = |\{i : \epsilon_i=-1\}|$. We obtain a presentation
$$x=\delta^{-n} b_1\cdots b_m,$$
with each $b_i$ being (up to an application of $\tau$ an integer number of times, which
preserves length by Proposition \ref{samelen})
$a_i$ if $\epsilon_i=1$, and $\bar a_i$ otherwise.
In particular, $\lmin(b_i)=1$ if $\epsilon_i=1$, and $\lmin(\bar a_i)=\lmin(\delta)-1$ otherwise.

Let $\delta^k s_1\cdots s_j$ be the left-weighted form of $b_1\cdots b_m$.
Then the greedy normal form of $x$ is
$\delta^{-n+k} s_1\cdots s_j$,
which cannot be longer than $\delta^{-n}\delta^k \allowbreak s_1\cdots s_j$.
As expressions of positive elements all have the same length,
the length of $\delta^k \allowbreak s_1\cdots s_j$ is exactly that of
$b_1\cdots b_m$.
Thus,
\begin{eqnarray*}
\lG(x) & \le & n\lmin(\delta)+\lmin(b_1\cdots b_m) = n\lmin(\delta)+\lmin(b_1\cdots b_m) =\\
 & = & n\lmin(\delta)+n(\lmin(\delta)-1)+(m-n) = \\
 & = & n(2\lmin(\delta)-2)+m\le (2\lmin(\delta)-1)m,
\end{eqnarray*}
as $n\le m$.\footnote{The step before last is added to emphasize that for random words, the
upper bound is far from being optimal. Indeed, in this case we have $n\approx m/2$, which gives roughly
half of the mentioned bound. There is an elbow room for improvements in the random case.}

(4) This can be proved as in the proof of (3).
Alternatively, one can use Charney's Theorem \cite{Char95}, extended to general Garside groups
by Dehornoy and Paris \cite{DP}, that the number of simple elements in the rational
normal form is minimal amongst presentations of $x$ as a product of simple elements (possibly inverted):
If $x\in M^{\pm 1}$, we can use (1) or (2) and there is nothing to prove.
Otherwise, let $x=a_1^{\epsilon_1}\cdots a_m^{\epsilon_m}$ be a minimal presentation of $x$.
In particular each $a_i^{\epsilon_1}$ is a (possibly inversed) simple element.
Thus, the number $n$ of simple elements in the rational form of $x$ is at most $m$.
As $x\notin M^{\pm 1}$, no simple element in the rational form of $x$ is $\delta$.
It follows that $\lR(x)\le (\lmin(\delta)-1)m$.

\medskip

(1) shows that the lower bounds cannot be improved.
To see that the upper bounds in (3) cannot be improved, consider
$\lG(a^{-m})$ for $m$ positive and an atom $a$
\end{proof}

The following corollary of Theorem \ref{basiclengths} is of special interest.
In 1994, Berger supplied an efficient method to compute a minimal length representative
of an element of $B_3$, in terms of Artin generators \cite{Berger}.
We show that the same is true for the BKL presentation. Indeed, a minimal
length representative for the BKL presentation is supplied by the rational normal form.

\bcor
Consider the BKL presentation of $B_3$. For each $x\in B_3$,
$\lR(x)=\lmin(x)$.
\ecor
\bpf
Here, $\lmin(\delta)=2$. By Theorem \ref{basiclengths},
$\lmin(x)\le \lR(x)\le (\lmin(\delta)-1)\lmin(x) = \lmin(x)$.
\epf

\brem\label{support}
Let $M$ be a Garside monoid, and $G$ be its fractions group.
Dehornoy and Paris \cite{DP} proved that for each $x\in G$, there is
a unique pair $(u,v)\in M^2$, such that $x=u\inv v$.
It follows that for each braid $x$, the rational normal form of $x$ belongs
to $B_N$ with the smallest possible $N$. In fact, if we define the \emph{support} of
a braid as the set of strands that cross in every braid representative, then
the rational normal form of $x$ detects its support. This is another reason why
rational normal forms approximate the minimal length.
\erem

\brem
We do not know whether the upper bound in (4) of Theorem \ref{basiclengths} can be improved.
At first it seems that for positive $m$ and distinct non-commuting atoms $a,b$,
$\lR(a^mb^{-m})=(\lmin(\delta)-1)\lmin(a^mb^{-m})$, but this is not the case:
Consider $\sigma_2^{2}\sigma_1^{-2}$ in the Artin presentation of $B_N$. Its rational normal form in $B_3$
(and thus by Remark \ref{support} in $B_N$ for all $N$) is
$(\sigma_1\inv \sigma_2\inv)\cdot(\sigma_2\inv\sigma_1\inv)\cdot(\sigma_2\sigma_1)\cdot(\sigma_1\sigma_2)$,
and thus $\lR(x) = 8 = 2 · \lmin(x)$. But $\lmin(\Delta)-1=2$ only when $N=3$.
\erem

Theorem \ref{basiclengths} shows that $\lR$ gives a better approximation than $\lG$,
and gives a theoretical motivation for the results described in \cite{LBCS}.
Having both experimental \cite{LBCS} and theoretical evidence for the superiority of $\lR$ over
$\lG$, we concentrate henceforth on the former.

\subsection{Quasi-geodesics in embedded Garside groups}
We need not stop here, and may consider, as in the case of $B_N$, two
distinct Garside structures of the same group, such that one of them embeds in
the other. Let $M_1,M_2$ be Garside monoids with fundamental elements $\Delta,\delta$,
respectively, such that each atom of $M_1$ is also an atom of $M_2$, and the group
of fractions of $M_1$ coincides with that of $M_2$.
Then we may take a length in one Garside structure as an estimation for the length in the other.
We will denote the used structure by a superscripted index.
By Theorem \ref{basiclengths},
\begin{eqnarray*}
\lR^2(x) & \le & (\lmin^2(\delta)-1)\lmin^2(x)\le (\lmin^2(\delta)-1)\lmin^1(x);\\
\lR^1(x) & \le & (\lmin^1(\Delta)-1)\lmin^1(x).
\end{eqnarray*}
Thus, if $\lmin^2(\delta)<\lmin^1(\Delta)$, $\lR^2(x)$ has a smaller approximation factor at its upper bound.

For the lower bound, let $A_2$ be the set of atoms of $M_2$, and set
$$\alpha=\max\{\lmin^1(a) : a\in A_2\}.$$
Then $\lmin^1(x)\le\alpha\lmin^2(x)$, and thus
$$\lmin^1(x)\le\alpha\lmin^2(x)\le\alpha\lR^2(x).$$
This gives the following.
\bthm\label{sym}
In the above notation,
$$\frac1\alpha\lmin^1(x)\le\lR^2(x)\le(\lmin^2(\delta)-1)\lmin^1(x).\qed$$
\ethm

The advantage of Theorem \ref{sym} is that the distortion factors are
symmetrized around the used length function $\lR^2(x)$.
Our main application is the following.

\subsection{The case of the braid group}

Consider the braid group as generated by the Artin monoid $B_N^+$
as well as by the BKL monoid $BKL_N^+$ (Examples \ref{ArtinBN}--\ref{BKLBN}),
and let $\Delta$ and $\delta$ be their respective fundamental elements.
Consider the minimal lengths $\lmin^1$ for the Artin structure, and $\lmin^2$ for
the BKL structure of $B_N$, respectively.

$\lmin^1(\Delta)=N(N-1)/2$, whereas $\lmin^2(\delta)=N-1$.
For each atom $a_{t,s}$ of $BKL_N^+$, $\lmin^1(a_{t,s})\le 2(t-s-1)+1 = 2(t-s)-1$.
In particular, the maximum $\alpha$ of all these lengths satisfies
$$\alpha\le 2N-3.$$
By Theorem \ref{sym}, we have that $\lR^2$, the length in BKL generators of
the rational normal form in the BKL structure of $B_N$, is quite symmetrically
close to the minimal \emph{Artin} length:

\bcor\label{app}
For each $x\in B_N$:
$$\frac1{2N-3}\lmin^1(x)\le\lR^2(x)\le(N-2)\lmin^1(x).\qed$$
\ecor

For comparison, measuring the minimal Artin length by working solely with the
Artin structure of $B_N$, we only have (by Theorem \ref{basiclengths}):
$$\lmin^1(x)\le\lR^1(x)\le (\lmin^1(\Delta)-1)\lmin^1(x) = \frac{N^2-N-2}{2}\lmin^1(x).$$
The gain may be viewed as follows:
In the latter case, we have a constant (in $N$) error factor from below, and quadratic
error from above. In Corollary \ref{app}, both errors are linear, that is, the errors
are symmetrized by dividing by $O(N)$ terms.

Another matter, which we cannot prove at present, is that the lower bound in
Corollary \ref{app} seems to be a big underestimate in the generic case.
It seems to us that in the generic case, the lower bound factor
should not be much smaller than 1 (indeed, it may be greater than 1).

In summary, we have theoretical evidence suggesting that
estimating the minimal length in Artin generators by using rational BKL normal form
should be better than the same estimation using rational Artin normal form.
We now turn to experimental results concerning the random case.

\section{Experimental results}\label{compar}

\subsection{Initial experiments}
For the Artin presentation, it is shown in \cite{LBCS} that
the rational Artin length is much better than greedy Artin length, at least with
regards to solving random equations with difficult parameters.
Our initial experiments showed that this is also the case
for the BKL presentation:
The rational BKL length is better than greedy BKL length.

In the initial phase of this project, we have compared various length
functions induced by various alternative ways of measuring lengths of elements,
and found out that only the rational BKL length outperforms the rational Artin length when
the problem's parameters are getting difficult.
The remainder of this report is therefore dedicated to the comparison of the these two leading
candidates.

\subsection{A detailed comparison}
We adopt the basic framework of \cite{Anshel, BraidEqns, LBCS}:
The equations are in a finitely generated group
$G=\<a_1,\dots,a_\ng\>\le B_\ns$, where $\ns$ denotes
the number of strings and $\ng$ denotes the number of generators of $G$.
Each generator $a_i$ is a word in $B_\ns$ obtained by multiplying $\wl$
(\emph{word length}) independent uniformly random elements of $\{\sigma_1,\dots,\sigma_{\ns-1}\}^{\pm 1}$.
In $G$, we build a sentence $X$ of length $\sl$ (\emph{sentence length}):
$$X=a_1a_2\cdots a_\sl$$
(For the while, we restrict $\sl\le\ng$).
Some of the $a_i$-s may be equal, but we did not force
that intentionally.

We begin with a description of a test suitable for
groups $G$ which are close to being free.
For each $i\in\{1,\dots,\ng\}$ and each $\epsilon\in\{1,-1\}$,
we give the generator $a_i^{\epsilon}$ the score
$$\ell(a_i^{-\epsilon}X),$$
sort the generators according to their scores (position $1$ is for the
shortest length), and reorder each block of identical scores by applying a random permutation.
We then keep in a histogram the position of $a_1$.
We do one such computation for each sample of $G$ and $X$.

While $a_1a_2\cdots a_\sl$ is not the
way a random $\sl$ sentence in $G$ was defined, this
does not make the problem easier: We use each group $G$ to produce
only one such sentence.

To partially compensate for the fact that $G$ need not be free,
we do the following.
There could be several $i\in\{1,\dots,\ng\}$
such that $X=a_ia_1\cdots a_{i-1}a_{i+1}\cdots a_\sl$.
Let $\COR$ denote the set of these $a_i$, the
\emph{correct first generators}.
After sorting all generators as above, instead of looking
for the position of $a_1$, we look at the lowest position
an element of $\COR$ attained.

\begin{rem}
A more precise, but infeasible, way to construct $\COR$ would be
to find all \emph{shortest presentations} of $X$ as a product of elements from $\{a_1,\dots,a_m\}^{\pm 1}$,
and let $\COR$ be the set of the first generators in
these presentations. For the parameters we have checked,
we believe that this should not make a big difference.
The results in Section \ref{nssec} support this hypothesis.
\end{rem}

We have also checked one set of cases where $\sl>\ng$.
In these cases we defined
$$X=a_{i_1}a_{i_2}\cdots a_{i_\sl},$$
where $i_j= (j-1 \bmod \ng)+1$ for $j=1,\dots,\sl$,
and made the obvious adjustments.

In summary, for each set of parameters $(\ns,\wl,\ng,\sl)$ mentioned below,
and for $\ell$ being either the rational Artin or the rational BKL length,
we have repeated the following at least $1,000$ times:
Choose $a_1,\dots,a_\ng$, compute $X$, compute $\COR$,
sort all generators $a_i^{\epsilon}$ according to the lengths $\ell(a_i^{-\epsilon}X)$,
find the lowest position attained by an element of $\COR$, and store this position number
in the histogram.

After dividing the numbers in the histogram by the numbers of samples made,
we obtain the distribution of the best position of a correct generator.
In light of the intended application described in the first two sections,
a natural measure to the effectiveness of $\ell$ is the graph of the
accumulated probability, showing for each $x=1,\dots,2\ng$ the probability
that some correct generator attained a position $\le x$.

The results of our experiments are divided into $4$ sets
such that in each set of experiments, only one parameter varies.
This shows the effect of that parameter on the difficulty
of the problem.
The varying parameter takes $3$ possible values, so we have $3$ \emph{pairs}
(since there are two length functions) of graphs. Each pair of graphs has its own
line style, so to allow plotting all $6$ graphs on the same figure.

For all pairs, one of the graphs is always above or almost the same as the other.
Fortunately, in all cases, \emph{it is the rational BKL length which is above the
rational Artin length},
so there is no need to supply this information in the figure.

Finally, since the accumulated distributions all reach $1$ for $x=2\ng$,
the graphs are more interesting for the smaller values of $x$.
We therefore plot only the first $35$ values of $x$.

\subsection{When the sentence length varies}
Fix $\ns=64, \wl=8, \ng=128$. Figure \ref{sl} shows the accumulated probabilities for
$\sl\in\{32,64,128\}$.
\myfigure{
\begin{changemargin}{-2cm}{-2cm}
\begin{center}
\Graphics{height=8cm, width=12.5cm}{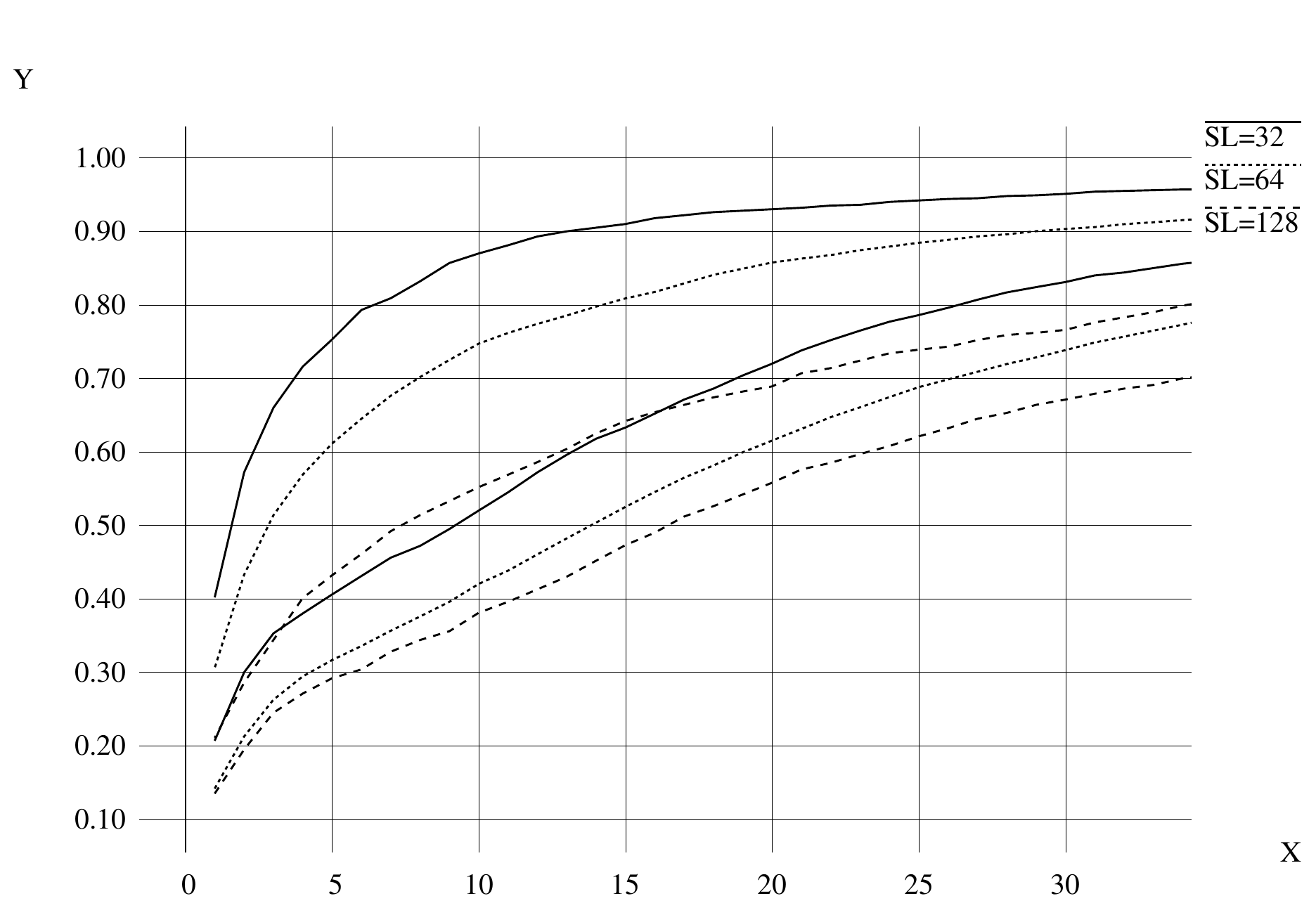}
\caption{When $\sl$ varies}\label{sl}
\end{center}
\end{changemargin}
}

\subsection{When the word length varies}
For $\ns=\sl=64,\ng=128$, and $\wl\in\{8,16,32\}$,
we obtain the graphs in Figure \ref{wl}.
The problem gets easier when $\wl$ increases, since
this way $G$ gets closer to a free group (where the length
approach is optimal). The remarkable observation is that
the harder the problem becomes (by making $\wl$ smaller),
the greater the improvement of the rational BKL length over
the rational Artin length becomes.
\myfigure{
\begin{changemargin}{-2cm}{-2cm}
\begin{center}
\Graphics{height=8cm, width=12.5cm}{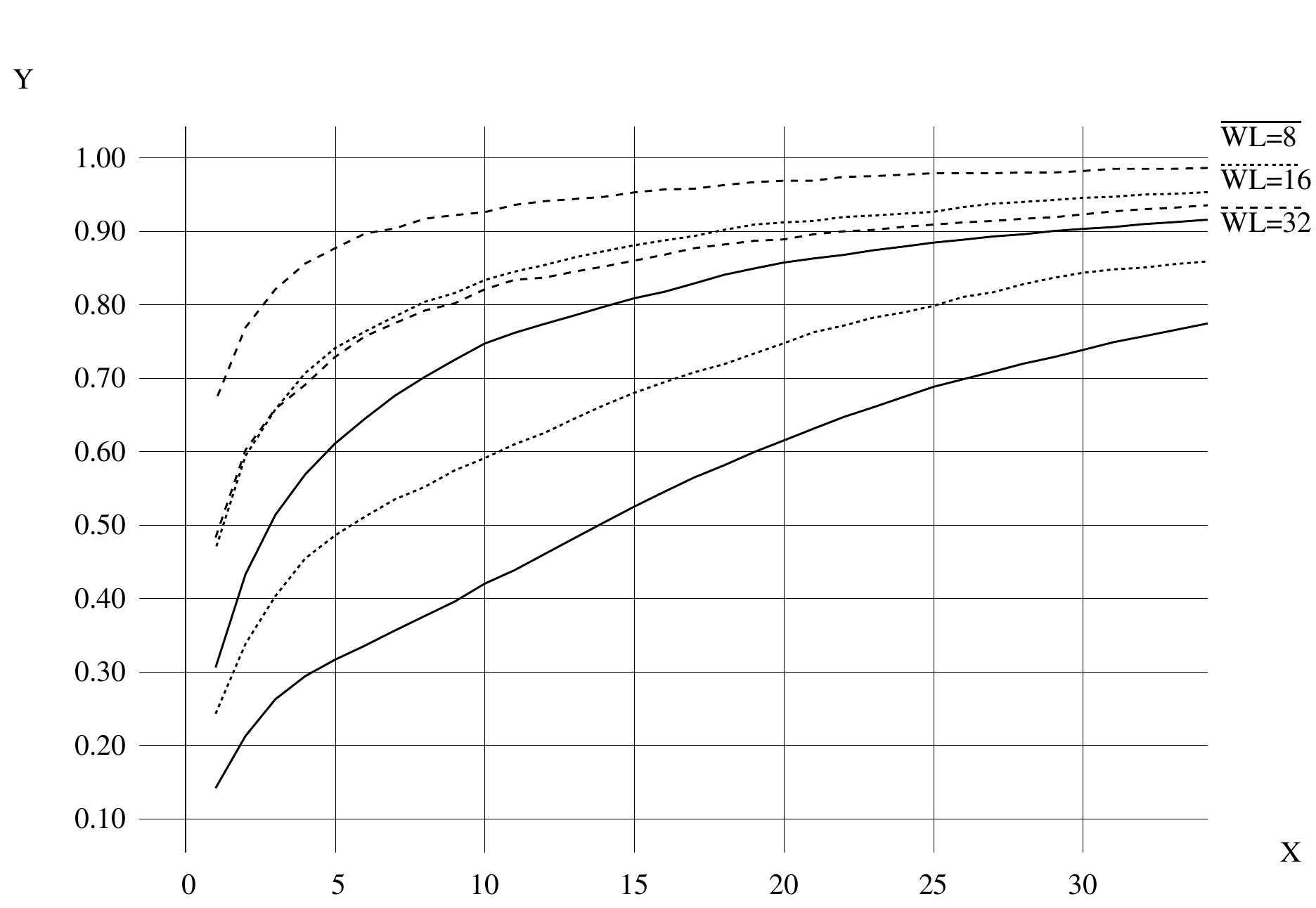}
\caption{When $\wl$ varies}\label{wl}
\end{center}
\end{changemargin}
}

\subsection{When the number of generators varies}
Now set  $\ns=\sl=64,\wl=8$, and let $\ng\in\{32,64,128\}$.
The graphs appear in Figure \ref{ng}.
Here too, the more difficult the problem becomes (by increasing
the number of generators), the greater the advantage of BKL over
Artin is. Moreover, the graphs show that doubling $\ng$ has
little influence on the performance of the rational BKL length,
whereas it seriously degrades the performance of the rational Artin
length.
\myfigure{
\begin{changemargin}{-2cm}{-2cm}
\begin{center}
\Graphics{height=8cm, width=12.5cm}{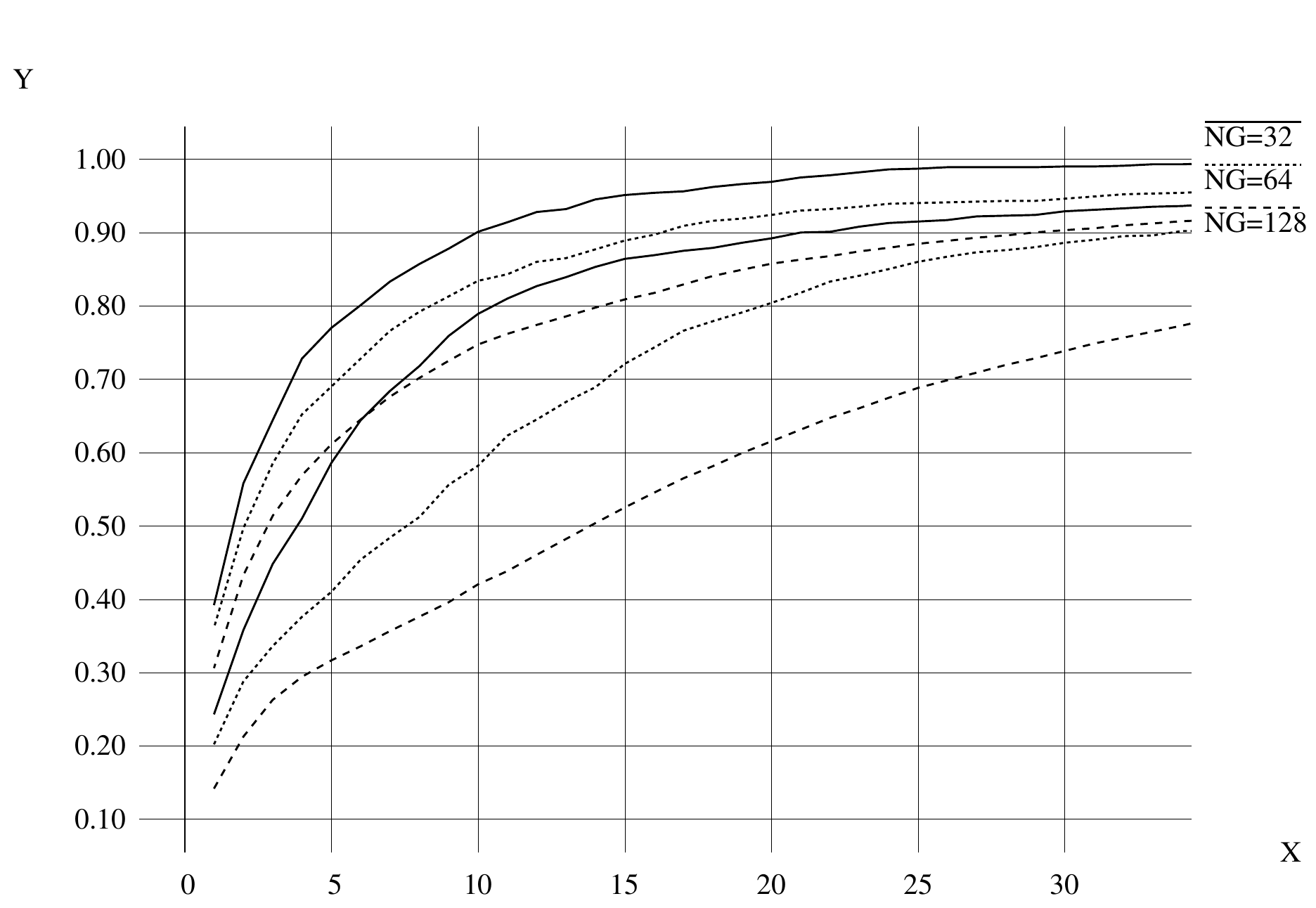}
\caption{When $\ng$ varies}\label{ng}
\end{center}
\end{changemargin}
}

\subsection{When the number of strings varies}\label{nssec}
Finally, set  $\wl=8,\sl=64,\ng=128$, and let $\ns\in\{16,32,64\}$.
Here, the problem becomes easier when we increase $\ns$ (Figure \ref{ns}).
This is \emph{not} in accordance with earlier results in \cite{LBCS, BraidEqns},
and is perhaps due to the fact that we allow any correct generator, whereas in
the earlier works we only counted $a_1$ a success.
Indeed, the more strings there are, the greater the chances are that
words of length $8$ commute. On the other hand, the graphs show that
while the BKL approach benefits a great deal when the number of strings is doubled,
this is not quite so for the Artin approach.
This means that the improvement in success rates due to commuting generators
is not substantial.

\myfigure{
\begin{changemargin}{-2cm}{-2cm}
\begin{center}
\Graphics{height=8cm, width=12.5cm}{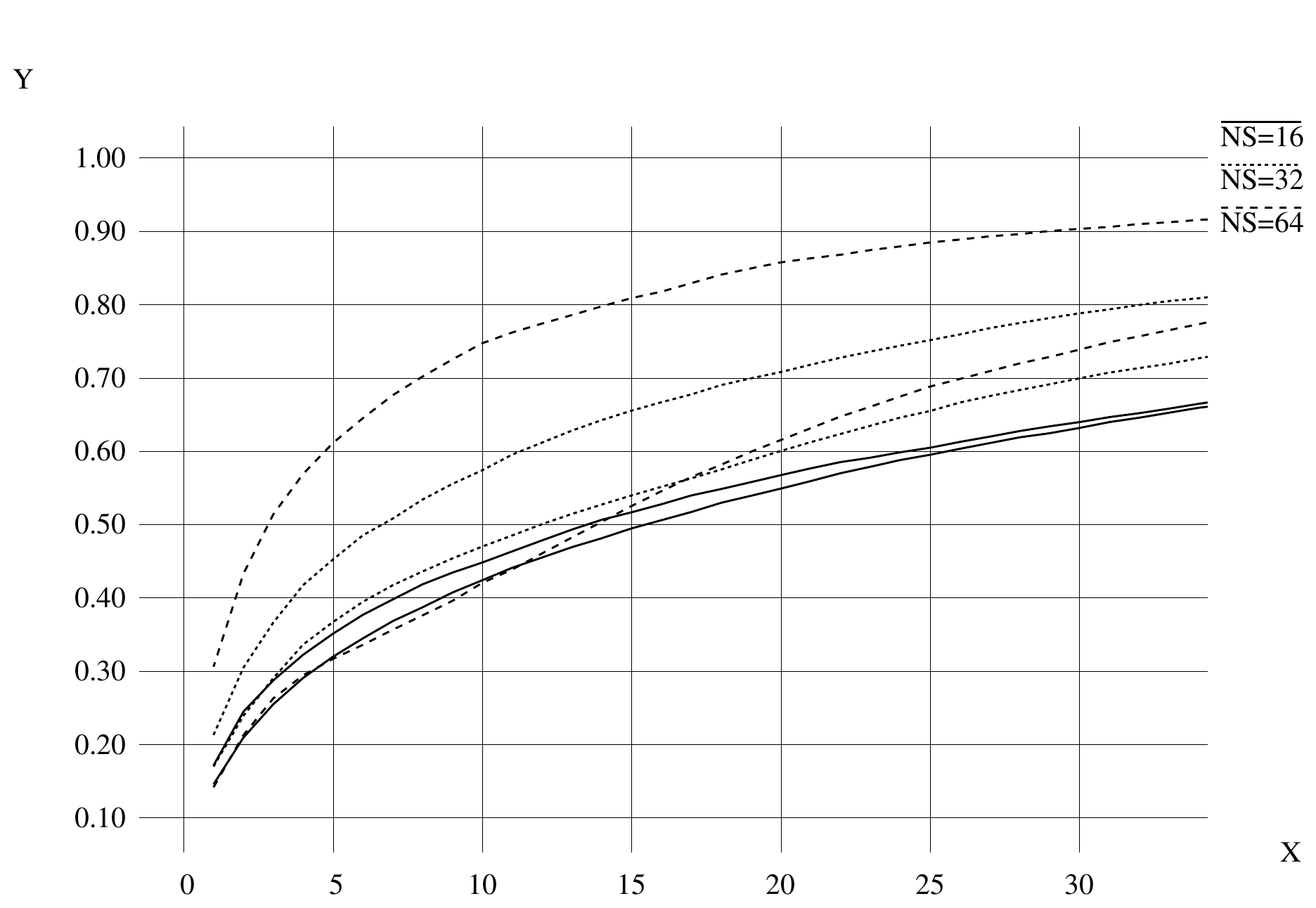}
\caption{When $\ns$ varies}\label{ns}
\end{center}
\end{changemargin}
}

\section{Concluding remarks and proposed future research}

Memory-length algorithms give a powerful heuristic method to solve
\emph{arbitrary} equations in noncommutative groups, and consequently a variety of
otherwise intractable problems.
These algorithms rely on a good length function on the group in question.
In the past, \emph{greedy Artin length} was used as a length function on the braid group,
and it was realized that \emph{rational Artin length} gives better results.

In this paper, we suggested to use \emph{rational BKL length} to measure the minimal
Artin length, and gave theoretical as well as experimental evidence for the advantage
of the new function over rational Artin length,
at least when randomization is modelled as in \cite{Anshel}.

The main drawback in our estimations is that they give much larger lengths than
the minimal length. Some interesting directions for possible improvements are:
\be
\itm As we have seen, the rational form can be computed from the greedy normal from
by ``removing'' $\delta$-s from the leading simple elements. We may be more greedy,
and remove the available $\delta$-s from the (leftmost) \emph{longest} simple elements in
the greedy normal form.\footnote{This was suggested to us by Uzi Vishne.}
This gives a new normal form in $B_N$, which has shorter length in terms
of atoms. The resulting length function may be yet better than the one proposed here.
\itm For each $x$ and each proposal for a length function of $x$, we can take the minimum
of the lengths of several elements whose minimal length is not smaller than that of $x$,
including:
$x$, $x\inv$, $x^{\delta^k}$ for each $k=1,\dots,m-1$, where
$m$ is the minimal with $\delta^m$ central.
\itm Since we use left-oriented normal forms in our estimations, we can also
try the corresponding \emph{right-oriented} normal forms, and take the minimum.
\itm We can iterate conjugation by $\delta$ and inverses (and other operations which are not increasing the
minimal length) with shortening heuristics like Dehornoy handle-reduction.
In \cite{MSU05} this was done only to a very limited extent.
\itm In \cite{MSU05}, Dehornoy handle-reduction was applied to the greedy normal form to obtain an estimation
of the minimal length. We conjecture that applying Dehornoy handle-reduction to the rational normal
form would give better estimations.
\ee

\subsection*{Acknowledgements}
We thank Joan Birman and Dima Ruinskiy for their comments on earlier versions of
the paper.
We also thank Patrick Dehornoy and Sang Jin Lee for
informative discussions concerning our notation, and the referees for their useful
comments.
A special thanks is owed to Arkadius Kalka for useful discussions and suggestions.

\end{document}